\documentclass[fleqn]{mat01}
\usepackage{times,mathtimy,amssymb,latexsym}
\begin{document}

\setcounter{page}{371} \firstpage{371}

\newcommand{\RR}{{\mathbb R}}
\newcommand{\DD}{{\mathbb D}}
\newcommand{\CC}{{\mathbb C}}
\newcommand{\ZZ}{{\mathbb Z}}
\newcommand{\QQ}{{\mathbb Q}}
\newcommand{\NN}{{\mathbb N}}
\def\a{\alpha}
\def\vp{\varphi}
\def\ol{\overline}
\def\a{\alpha}
\def\om{\omega}
\def\b{\beta}
\def\g{\gamma}
\def\ol{\overline}
\def\ve{\varepsilon}
\def\de{\delta}
\def\pt{\partial}
\def\qand{\quad\mbox{ and }\quad}
\def\qfor{\quad\mbox{ for }\quad}
\def\d{\mbox{\rm d}}
\def\e{\mbox{\rm e}}

\newtheorem{theore}{Theorem}
\renewcommand\thetheore{\arabic{section}.\arabic{theore}}
\newtheorem{theor}[theore]{\bf Theorem}
\newtheorem{lem}[theore]{\it Lemma}
\newtheorem{coro}[theore]{\rm COROLLARY}

\renewcommand{\thefootnote}{\fnsymbol{footnote}} 

\title{Weighted composition operators from Bergman-type spaces into Bloch spaces}

\markboth{Songxiao Li and Stevo Stevi\'c}{Weighted composition operators from
Bergman-type spaces}

\author{SONGXIAO LI$^{1,2}$ and STEVO STEVI\'C$^{3}$}

\address{$^{1}$Department of Mathematics, Shantou University, 515~063 Shantou, Guangdong, China\\
\noindent $^{2}$Department of Mathematics, Jia Ying University, 514~015 Meizhou, Guangdong, China\\
\noindent $^{3}$Mathematical Institute of the Serbian Academy of Science,
Knez Mihailova 35/I, 11000 Beograd, Serbia\\
\noindent E-mail: jyulsx@163.com; lsx@mail.zjxu.edu.cn; sstevic@ptt.yu;
sstevo@matf.bg.ac.yu}

\volume{117}

\mon{August}

\parts{3}

\pubyear{2007}

\Date{MS received 5 June 2006; revised 3 October 2006}

\begin{abstract}
Let $\varphi$ be an analytic self-map and $u$ be a fixed analytic function on the open
unit disk $D$ in the complex plane $\CC.$ The weighted composition operator is defined\break by
\begin{equation*}
uC_\varphi f =u \cdot (f\circ \varphi), \ \ f \in H(D).
\end{equation*}
Weighted composition operators from Bergman-type spaces into
Bloch spaces and little Bloch spaces are characterized by function
theoretic properties of their inducing maps.
\end{abstract}


\maketitle

\section{Introduction}\vspace{.5pc}

Let $D$ be the open unit disk in the complex plane $\CC.$ Denote
by $H(D)$ the class of all functions analytic on $D$. An analytic
self-map $\varphi\hbox{\rm :}\ D\rightarrow D$ induces the
composition operator $C_\varphi$ on $H(D)$, defined by
$C_\varphi(f) = f( \varphi(z))$ for $f$ analytic on $D$. It is a
well-known consequence of Littlewood's subordination principle
that the composition operator $C_\varphi$ is bounded on the
classical Hardy and Bergman spaces (see, for example~\cite{cm}).

Recall that a linear operator is said to be bounded if the image
of a bounded set is a bounded set, while a linear operator is
compact if it takes bounded sets to sets with compact closure. It
is interesting to provide a function theoretic characterization of
when $\varphi$ induces a bounded or compact composition operator
on various spaces. The book \cite{cm} contains plenty of
information on this topic.

Let $u$ be a fixed analytic function on the open unit disk. Define
a linear operator $uC_\varphi$ on the space of analytic functions
on $D$, called a weighted composition operator, by $uC_\varphi f
=u \cdot (f\circ \varphi),$ where $f$ is an analytic function on
$D$. We can regard this operator as a generalization of a
multiplication operator and a composition operator.

A positive continuous function $\phi$ on $[0,1)$ is called normal,
if there exist positive numbers $s$ and $t,$ $0<s<t,$ such that
\begin{equation*}
\frac{\phi(r)}{(1-r)^s}\downarrow 0, \qquad
\frac{\phi(r)}{(1-r)^t} \uparrow\infty
\end{equation*}
as $r\rightarrow 1^-$ (see, for example~\cite{hu,sw}).\pagebreak

For $0<p<\infty$, $0<q<\infty$ and a normal function $\phi$, let
$H(p,q,\phi)$ denote the space of all analytic functions $f$ on
the unit disk $D$ such that
\begin{equation*}
\|f\|_{p,q,\phi}=\left(\int_0^1
M_q^p(r,f)\frac{\phi^{p}(r)}{1-r}r\d r\right)^{1/p}<\infty,
\end{equation*}
where the integral means $M_p(f,r)$ are defined by
\begin{equation*}
M_p(f,r)=\left(\frac{1}{2\pi} \int_0^{2\pi} |f(r\e^{i\theta})|^p
\d\theta\right)^{1/p},\qquad 0\leq r<1.
\end{equation*}
For $1\leq p<\infty$, $H(p,q,\phi)$, equipped with the norm
$\|\cdot\|_{p,q,\phi}$ is a Banach space. When $0<p<1$, $
\|f\|_{p,q,\phi}$ is a quasinorm on $H(p,q,\phi)$, $H(p,q,\phi)$
is a Frechet space but not a Banach space. If $0<p=q<\infty$, then
$H(p,p,\phi)$ is the Bergman-type space
\begin{equation*}
H(p,p,\phi)=\left\{f \in H(D)\hbox{\rm :}\  \int_D |f(z)|^p
\frac{\phi^p(|z|)}{1-|z|} \d A(z)<\infty\right\}.
\end{equation*}
Here $\d A$ denotes the normalized Lebesgue area measure on the unit
disk $D$ such that $A(D) = 1$. Note that if $\phi(r)=(1-r)^{1/p}$,
then $H(p,p,\phi)$ is the Bergman space $A^p$.

An analytic function $f$ in $D$ is said to belong to the Bloch
space $\mathcal{B}$ if
\begin{equation*}
B(f) = \sup_{z \in
D}(1-|z|^2)|f'(z)|<\infty.
\end{equation*}
The expression $B(f)$ defines a seminorm while the natural norm is
given by $\|f\|_{\mathcal{B}}=|f(0)|+B(f) $. The norm makes
$\mathcal{B}$ into a conformally invariant Banach space. Let
$\mathcal{B}_0$ denote the subspace of $\mathcal{B}$ consisting of
those $f \in \mathcal{B}$ for which $(1-|z|^2)|f'(z)|\rightarrow 0 $, as $|z| \rightarrow 1$. This
space is called a little Bloch space. For more information on
Bloch spaces see, for example~\cite{cm,r,ss,zhu,zhu1} and
the references therein.

In \cite{o}, Ohno has characterized the boundedness and
compactness of weighted composition operators between $H^\infty$,
the Bloch space $\mathcal{B}$ and the little Bloch space
$\mathcal{B}_0$. In \cite{oz}, Ohno and Zhao have characterized
the boundedness and compactness of weighted composition operators
on the Bloch space. Weighted composition operators between
Bloch-type spaces are characterized in \cite{osz} (see also
\cite{mr}). In the setting of the unit ball or the unit polydisk,
some necessary and sufficient conditions for a composition
operator or weighted composition operator to be bounded or compact
are given, for example, in \cite{cm,sl,zs2,zhou}.

In this paper we study the weighted composition operators from the
Bergman-type space $H(p,p,\phi)$ into the Bloch space
$\mathcal{B}$ and the little Bloch space $\mathcal{B}_0$. As
corollaries, we obtain the complete characterizations of the
boundedness and compactness of composition operators from Bergman
spaces into Bloch spaces.

In this paper, constants are denoted by $C$, they are positive and
may differ from one occurrence to the next. The notation $a \preceq
b$ means that there is a positive constant $C$ such that $a \leq C
b$. If both $a\preceq b$ and $b\preceq a$ hold, then one says that
$a \asymp b$.

\section{Auxiliary results}

In this section, we give some auxiliary results which will be used in proving
the main results of the paper. They are incorporated in the lemmas which follow.

\begin{lem} Let $0<p<\infty$. If $f \in
H(p,p,\phi),$ then
\begin{equation}
|f(z)| \leq C
\frac{ \|f\|_{H(p,p,\phi)}} {\phi(|z|)(1-|z|^2)^{1/p}}, \qquad z
\in D.
\end{equation}
\end{lem}

\begin{proof}
Let $\beta(z,w)$ denote the Bergman metric between
two points $z$ and $w$ in $D$. It is well-known that
\begin{equation*}
\beta(z,w)=\frac{1}{2}\log\frac{1+|\varphi_z(w)|}{1-|\varphi_z(w)|}.
\end{equation*}
For $a\in D$ and $r>0$ the set $D(a,r)=\{z\in D\hbox{\rm :}\ \beta(a,z)<r\}$ is
the Bergman metric disk centered at $a$ with radius $r$. It is
well-known that (see \cite{zhu})
\begin{equation} \frac{(1-|a|^2)^2}{|1-  \bar{a}z|^4}\asymp
\frac{1}{(1-|z|^2)^2} \asymp \frac{1}{(1-|a|^2)^2} \asymp
\frac{1}{|D(a,r)|},
\end{equation}
when $z\in D(a,r)$, where $|D(a,r)|$ denotes the area of the disk
$D(a,r).$

From (2) and since $\phi(r)$ is normal it is not difficult to see
that for a fixed $r\in(0,1)$ the following relationship holds:
\begin{equation}
\phi(|z|)\asymp \phi(|a|),\quad z\in D(a,r).
\end{equation}

For $0<r<1$ and $z \in D$, by the subharmonicity of $|f(z)|^p$,
(2) and (3), we have that
\begin{align*}
|f(z)|^p&\leq \frac{C}{(1-|z|^2)^2}\int_{D(z,r)}
|f(w)|^p\d A(w)\\[.5pc]
&\leq  \frac{C}{(1-|z|^2) \phi^p(|z|)}\int_{D(z,r)}
|f(w)|^p\frac{\phi^p(|w|)}{1-|w|}\d A(w)\\[.5pc]
& \leq \frac{C}{(1-|z|^2) \phi^p(|z|)}\int_D
|f(w)|^p \frac{\phi^p(|w|)}{1-|w|}\d A(w)\\[.5pc]
& \leq \frac{C\|f\|^p_{H(p,p,\phi)}}{(1-|z|^2) \phi^p(|z|)},
\end{align*}
from which the desired result follows.
\end{proof}

The following lemma can be found in \cite{hu}.

\begin{lem} Let $0<p<\infty$. Then for $f \in
H(D),$
\begin{equation*}
\|f\|^p_{H(p,p,\phi)}\asymp |f(0)|^p+\int_D |f'(z)|^p(1-|z|^2)^p
\frac{\phi^p(|z|)}{1-|z|}\d A(z).
\end{equation*}
\end{lem}

\begin{lem} Let $0<p<\infty$. If $f \in
H(p,p,\phi)$ and $z\in D,$ then
\begin{equation} |f'(z)| \leq C
\frac{ \|f\|_{H(p,p,\phi)}} {\phi(|z|)(1-|z|^2)^{1/p+1}}, \qquad z
\in D.
\end{equation}
\end{lem}

\begin{proof}
By the subharmonicity of $|f'(z)|^p$, (2) and (3),
and Lemma~2.2 we have that
\begin{align*}
\hskip -4pc |f'(z)|^p&\leq \frac{C}{(1-|z|^2)^2}\int_{D(z,r)}
|f'(w)|^p\hbox{d}A(w)\\[.7pc]
\hskip -4pc &\leq \frac{C}{(1-|z|^2)^{p+1}\phi^p(|z|)}\int_{D(z,r)}
\frac{\phi^p(|w|)}{1-|w|}(1-|w|)^p|f'(w)|^p\hbox{d}A(w)\\[.7pc]
\hskip -4pc & \leq  \frac{C}{(1-|z|^2)^{p+1}\phi^p(|z|)}\int_D
\frac{\phi^p(|w|)}{1-|w|}(1-|w|)^p|f'(w)|^p\hbox{d}A(w)\\[.7pc]
\hskip -4pc &\leq \frac{C\|f\|^p_{H(p,p,\phi)}}{(1-|z|^2)^{p+1} \phi^p(|z|)},
\end{align*}
from which the result follows.
\end{proof}

The following lemma can be found in \cite{sw}.

\begin{lem}
For $\beta>-1$ and $m>1+\beta$ we
have
\begin{equation*}
\int_0^1 \frac{(1-r)^\beta}{(1-\rho r)^m}\d r\leq
C(1-\rho)^{1+\beta-m},\qquad 0<\rho<1.
\end{equation*}
\end{lem}

The following criterion for compactness follows by standard arguments similar, for
example, to those outlined in Proposition~3.11 of \cite{cm}.

\begin{lem} The operator $uC_\vp\hbox{\rm :}\ H(p,p,\phi)\to {\cal B}$ is
compact if and only if for any bounded sequence $(f_n)_{n\in {\NN}}$ in $H(p,p,\phi)$
which converges to zero uniformly on compact subsets of $D,$ we have $\|uC_\vp
f_n\|_{{\cal B}}\to 0$ as $n\to\infty.$
\end{lem}

\section{The boundedness and compactness of the operator $\pmb{uC_{\varphi}\hbox{\rm :}\ H(p,p,\phi)\rightarrow
\mathcal{B}}$}

In this section we characterize the boundedness and compactness of the weighted
composition operator $uC_{\varphi}\hbox{\rm :}\ H(p,p,\phi)\rightarrow \mathcal{B}$.

\setcounter{theore}{0}
\begin{theor}[\!]
Suppose that $\varphi$ is an analytic self-map of the unit disk{\rm ,} $u \in H(D),$
$0< p<\infty$ and that $\phi$ is normal on $[0,1)$. Then{\rm ,} $uC_{\varphi}\hbox{\rm
:}\ H(p,p,\phi) \rightarrow \mathcal{B}$ is bounded if and only if the following
conditions are satisfied$:$

\begin{enumerate}
\renewcommand\labelenumi{{\rm (\roman{enumi})}}
\leftskip .25pc
\item$\left.\right.$\vspace{-2pc}
\begin{equation}
\hskip -1.25pc \sup_{z\in D}\frac{(1-|z|^2)|u'(z)| }{ \phi(|\varphi(z)|)
(1-|\varphi(z)|^2)^{1/p} }<\infty;
\end{equation}\vspace{.5pc}

\item$\left.\right.$\vspace{-2pc}
\begin{equation}
\hskip -1.25pc \sup_{z\in D}\frac{ (1-|z|^2) |u(z)\varphi'(z)| }{
 \phi(|\varphi(z)|)(1-|\varphi(z)|^2)^{1+1/p} }<\infty.
\end{equation}
\end{enumerate}
\end{theor}

\begin{proof}
Suppose that the conditions (i) and (ii) hold.
For arbitrary $z$ in $D$ and $f \in H(p,p,\phi)$, by Lemmas~2.1 and 2.3 we have\pagebreak
\begin{align}
&\hskip -4pc (1-|z|^{2})|(uC_{\varphi}f)'(z)|\nonumber\\[.7pc]
&\hskip -4pc \quad\,\leq (1-|z|^{2})|u'(z)||f(\varphi(z))|+(1-|z|^{2})|f'(\varphi(z))||u(z)\varphi'(z)|\nonumber\\[.7pc]
&\hskip -4pc \quad\,\leq (1-|z|^{2})|u'(z)|\frac{ C \|f\|_{H(p,p,\phi) }
}{\phi(|\varphi(z)|)(1-|\varphi(z)|^{2})^{1/p}}\nonumber\\[.7pc]
&\hskip -4pc \qquad\,+C(1-|z|^2) |u(z)\varphi'(z)|  \frac{ \|f\|_{H(p,p,\phi)} }{
 \phi(|\varphi(z)|)(1-|\varphi(z)|^2)^{1+1/p}}\nonumber\\[.7pc]
&\hskip -4pc \quad\,\leq \left(  \frac{C(1-|z|^2)|u'(z)|}{\phi(|\varphi(z)|)
(1-|\varphi(z)|^2)^{1/p} }+ \frac{ C (1-|z|^2) |u(z)\varphi'(z)|
}{\phi(|\varphi(z)|)(1-|\varphi(z)|^2)^{1+1/p} } \right)
 \|f\|_{H(p,p,\phi)}.
\end{align}
Taking the supremum in (7) over $D$ and then using conditions (5) and (6) we obtain
that the operator $uC_{\varphi}\hbox{\rm :}\ H(p,p,\phi)\rightarrow \mathcal{B}$ is bounded.
\end{proof}

Conversely, suppose that $uC_{\varphi}\hbox{\rm :}\ H(p,p,\phi)\rightarrow
\mathcal{B}$ is bounded. Then, taking the functions $f(z)=z$ and
$f(z)=1$ we obtain that the quantities
\begin{equation*} \sup_{z\in
D}(1-|z|^{2})|u(z)\varphi'(z)+u'(z)\varphi(z)|\qand  \sup_{z\in
D}(1-|z|^{2})|u'(z)|
\end{equation*}
are finite. Using these facts and the boundedness of the function $\varphi(z),$ we
have that
\begin{equation}
\sup_{z\in D}(1-|z|^{2})|u(z)\varphi'(z)|<\infty.
\end{equation}
For fixed $w\in D$, take
\begin{equation}
f_w(z)= \frac{(1-|w|^{2})^{t+1}}{\phi(|w|) (1-\bar w z)^{1/p+t+1} }.
\end{equation}
By Lemma~1.4.10 of \cite{r}, we know that
\begin{equation*}
M_p(f_w,r)\leq C \frac{(1-|w|^2)^{t+1}}{\phi(|w|) (1-r|w|)^{t+1}}.
\end{equation*}
Since $\phi$ is normal, by Lemma~2.4, we obtain
\begin{align*}
\hskip -4pc \|f_w\|^p_{H(p,p,\phi)}&=\int_0^1 M_p^p(f_w,r)\frac{\phi^p(r)}{1-r}r\d r\\[.7pc]
\hskip -4pc &\leq C\int_0^1 \frac{(1-|w|^2)^{p(t+1)}}{
\phi^p(|w|)(1-r|w|)^{p(t+1)}}\frac{\phi^p(r)}{1-r}\d r\\[.7pc]
\hskip -4pc &\leq C\left(\int_0^{|w|} \frac{(1-|w|^2)^{p(t+1)}}{
\phi^p(|w|)(1-r|w|)^{p(t+1)}}\frac{\phi^p(r)}{1-r}\d r\right.\\[.7pc]
\hskip -4pc &\left.\quad\,+\int_{|w|}^1\frac{(1-|w|^2)^{p(t+1)}}{
\phi^p(|w|)(1-r|w|)^{p(t+1)}}\frac{\phi^p(r)}{1-r}\d r\right)
\end{align*}
\begin{align*}
\hskip -4pc \phantom{\|f_w\|^p_{p,p,\phi}}&\leq C\frac
{(1-|w|^2)^{p(t+1)}}{\phi^p(|w|)}\frac{\phi^p(|w|)}{(1-|w|^2)^{pt}}\int_0^{|w|}
\frac{(1-r)^{pt-1}}{(1-r|w|)^{p(t+1)}}\d r\\[.7pc]
\hskip -4pc &\quad\,+C\frac
{(1-|w|^2)^{p(t+1)}}{\phi^p(|w|)}\frac{\phi^p(|w|)}{(1-|w|^2)^{ps}}\int_{|w|}^1
\frac{(1-r)^{ps-1}}{ (1-r|w|)^{p(t+1)}}d\leq C.
\end{align*}
Therefore $f_w\in H(p,p,\phi),$ and moreover $\sup_{w\in D}\| f_w
\|_{H(p,p,\phi)}\leq C$. Hence, we have
\begin{align*}
\hskip -4pc C\|uC_{\varphi}\| &\geq \|f_{\vp(\lambda)}\|_{H(p,p,\phi)
}\|uC_{\varphi}\| \geq \|
uC_{\varphi}f_{\vp(\lambda)}\|_{\cal B} \\[.5pc]
\hskip -4pc &\geq \left|(1/p+t+1) \frac{(1-|\lambda|^{2}) |u(\lambda)
\overline{\varphi(\lambda)} \varphi'(\lambda)| }{ \phi(|\varphi(\lambda)|)
(1-|\varphi(\lambda)|^2)^{1+1/p}}- \frac{(1-|\lambda|^2)|u'(\lambda)|}
{\phi(|\varphi(\lambda)|) (1-|\varphi(\lambda)|^2)^{1/p}}\right|,
\end{align*}
for every $\lambda\in D,$ from which it follows that
\begin{align}
& \frac{(1-|\lambda|^2)|u'(\lambda)|} {\phi(|\varphi(\lambda)|)
(1-|\varphi(\lambda)|^2)^{1/p}}\nonumber\\[.7pc]
&\quad\,\leq C\| uC_{\varphi}\|+(1/p+t+1) \frac{(1-|\lambda|^{2}) |u(\lambda)
\overline{\varphi(\lambda)} \varphi'(\lambda)|}{ \phi(|\varphi(\lambda)|)
(1-|\varphi(\lambda)|^2)^{1+1/p}}.
\end{align}
Further, for $\lambda\in D$, take
\begin{align*}
g_\lambda(z)=\frac{(1-|\varphi(\lambda)|^{2})^{t+2}}{\phi(|\varphi(\lambda)|)
(1-\overline{\varphi(\lambda)} z)^{1/p+t+2} }
-\frac{(1-|\varphi(\lambda)|^{2})^{t+1}}{\phi(|\varphi(\lambda)|)
(1-\overline{\varphi(\lambda)} z)^{1/p+t+1} }.
\end{align*}
Then, $\sup_{\lambda\in D}\| g_\lambda \|_{H(p,p,\phi)}\leq C$,
$\;g_\lambda(\varphi(\lambda))=0$ and
\begin{equation*}
 g'_\lambda(\varphi(\lambda)) =
\frac{\overline{\varphi(\lambda)}}{\phi(|\varphi(\lambda)|)(1-|\varphi(\lambda)|^2)^{1+1/p}}.
\end{equation*}
Thus,
\begin{align*}
C\| uC_{\varphi}\|\geq \| uC_{\varphi}g_\lambda\|_{\mathcal{B}}
\geq  \frac{
(1-|\lambda|^2)|u(\lambda)\overline{\varphi(\lambda)}\varphi'(\lambda)
| }{ \phi(|\varphi(\lambda)|)(1-|\varphi(\lambda)|^2 )^{1+1/p} },
\end{align*}
i.e. we have
\begin{equation}
\sup_{ \lambda \in D}\frac{
(1-|\lambda|^2)|u(\lambda)\overline{\varphi(\lambda)}\varphi'(\lambda)
| }{ \phi(|\varphi(\lambda)|)(1-|\varphi(\lambda)|^2 )^{1+1/p}
}<\infty.
\end{equation}
Thus for a fixed $\delta$, $0<\delta<1$, by (11),
\begin{equation}
\sup_{ |\varphi(\lambda)|>\delta }   \frac{
(1-|\lambda|^2)|u(\lambda)|| \varphi'(\lambda) | }{
\phi(|\varphi(\lambda)|) (1-|\varphi(\lambda)|^2 )^{1+1/p} }
<\infty.
\end{equation}
For $\lambda\in D$ such that $ |\varphi(\lambda)|\leq \delta $,
since $\phi$ is normal, we have
\begin{equation}
\frac{ (1-|\lambda|^2)|u(\lambda)\varphi'(\lambda) | }{
\phi(|\varphi(\lambda)|)(1-|\varphi(\lambda)|^2 )^{1+1/p} } \leq
\frac{C}{(1-\delta^2)^{1+1/p}\phi(\delta)}
(1-|\lambda|^2)|u(\lambda)\varphi'(\lambda)|.
\end{equation}
Hence, from (8) and (13), we obtain
\begin{equation}
\sup_{ |\varphi(\lambda)|\leq \delta }\frac{
(1-|\lambda|^2)|u(\lambda)\varphi'(\lambda) | }{
\phi(|\varphi(\lambda)|)(1-|\varphi(\lambda)|^2 )^{1+1/p}
}<\infty.
\end{equation}
The inequality in (6) follows from (12) and (14). Taking the supremum in (10) over
$\lambda\in D$ and using (6), (5) follows. This completes the proof of the
theorem.\vspace{.2pc}

\begin{theor}[\!]
Suppose that $\varphi $ is an analytic self-map of the unit disk{\rm ,} $u \in H(D),$
$0< p<\infty,$ that $\phi$ is normal on $[0,1)$ and that $uC_{\varphi}\hbox{\rm :}\
H(p,p,\phi)\rightarrow \mathcal{B}$ is bounded. Then{\rm ,} $uC_{\varphi}\hbox{\rm :}\
H(p,p,\phi)\rightarrow \mathcal{B}$ is compact if and only if the following conditions
are\break satisfied$:$\vspace{.2pc}

\begin{enumerate}
\renewcommand\labelenumi{{\rm (\roman{enumi})}}
\leftskip .25pc
\item$\left.\right.$\vspace{-1.9pc}
\begin{equation}
\hskip -1.25pc \lim_{|\varphi(z)|\rightarrow 1}\frac{(1-|z|^2)|u'(z)| }{
\phi(|\varphi(z)|) (1-|\varphi(z)|^2)^{1/p} }   =0;
\end{equation}\vspace{.6pc}

\item$\left.\right.$\vspace{-1.9pc}
\begin{equation}
\hskip -1.25pc \lim_{|\varphi(z)|\rightarrow 1} \frac{ (1-|z|^2) |u(z)\varphi'(z)| }{
 \phi(|\varphi(z)|)(1-|\varphi(z)|^2)^{1+1/p} } =0.
\end{equation}
\end{enumerate}
\end{theor}

\begin{proof}
First assume that conditions (i) and (ii) hold. In order to prove that $uC_{\varphi}$
is compact, according to Lemma~2.5, it suffices to show that if $(f_{n})_{n\in\NN}$ is
a bounded sequence in $H(p,p,\phi)$ that converges to 0 uniformly on compact subsets
of $D$, then $\| uC_{\varphi}f_{n}\|_{\mathcal{B}}\rightarrow 0$ as $n\to\infty.$ Let
$(f_{n})_{n\in\NN}$ be a sequence in $H(p,p,\phi)$ with $\sup_{n\in\NN}\|
f_n\|_{H(p,p,\phi)}\leq L$ and suppose $f_n \rightarrow 0$ uniformly on compact
subsets of $D$ as $n\to\infty.$

By the assumptions of the theorem we have that for any
$\varepsilon>0$, there is a constant $\delta,\; 0<\delta<1$, such
that $\delta<|\varphi(z)|<1$ implies
\begin{equation*}
\frac{(1-|z|^2)|u'(z)| }{ \phi(|\varphi(z)|)
(1-|\varphi(z)|^2)^{1/p} } <\varepsilon/L
\end{equation*}
and
\begin{equation*}
\frac{ (1-|z|^2)
|u(z)\varphi'(z)| }{
 \phi(|\varphi(z)|)(1-|\varphi(z)|^2)^{1+1/p} }<\varepsilon/L.
\end{equation*}
Let $K=\{ w \in D\hbox{\rm :}\ |w|\leq \delta \}$. Note that $K$ is a compact subset
of $D.$ From this, since $\phi$ is normal and using estimates from Lemmas~2.1 and 2.3,
we have that
\begin{align*}
&\hskip -4pc \| uC_{\varphi}f_{n}\|_{\mathcal{B}}\\[.4pc]
&\hskip -4pc \quad\, = \sup_{z\in D}(1-|z|^{2})|(uC_{\varphi}f_{n})'(z)|+|u(0)f_n(\varphi(0))| \\[.65pc]
&\hskip -4pc \quad\,\leq\!\sup_{z\in
D}(1-|z|^{2})|u'(z)f_{n}(\varphi(z))|\!+\!\sup_{z\in
D}(1-|z|^{2})|u(z)f_n'(\varphi(z))\varphi'(z)|+|u(0)f_n(\varphi(0))|\\[.65pc]
&\hskip -4pc \quad\,\leq \sup_{\{z\in D\hbox{\rm :}\ \varphi(z)\in
K\}}(1-|z|^{2})|u'(z)f_{n}(\varphi(z))|+ \sup_{\{z\in
D\hbox{\rm :}\ \de\leq|\varphi(z)|<1\}}(1-|z|^{2})|u'(z)f_{n}(\varphi(z))|\\[.65pc]
&\hskip -4pc \qquad\, +\sup_{\{z\in D\hbox{\rm :}\ \varphi(z)\in
K\}}(1-|z|^{2})|u(z)\varphi'(z)||f_{n}'(\varphi(z))|
\end{align*}
\begin{align*}
&\hskip -4pc \qquad\, +\sup_{\{z\in D\hbox{\rm :}\
\de\leq|\varphi(z)|<1\}}(1-|z|^{2})|u(z)\varphi'(z)||f_{n}'(\varphi(z))|+|u(0)f_n(\varphi(0))|\\[.5pc]
&\hskip -4pc \quad\, \leq \|u\|_{\mathcal{B}} \sup_{w\in K}|f_{n}(w)|+C\sup_{\{z\in
D\hbox{\rm :}\ \de\leq|\varphi(z)|<1\}}\frac{(1-|z|^2)|u'(z)| }{ \phi(|\varphi(z)|)
(1-|\varphi(z)|^2)^{1/p} }\|f_n\|_{H(p,p,\phi)}\\[.5pc]
&\hskip -4pc \qquad\, +M\sup_{w\in K}|f'_{n}(w)|+C\sup_{\{z\in D\hbox{\rm :}\
\de\leq|\varphi(z)|<1\}}\frac{ (1-|z|^2) |u(z)\varphi'(z)| }{
 \phi(|\varphi(z)|)(1-|\varphi(z)|^2)^{1+1/p}}\|f_n\|_{H(p,p,\phi)}\\[.5pc]
&\hskip -4pc \qquad\, +|u(0)f_n(\varphi(0))|\\[.7pc]
&\hskip -4pc \quad\, \leq \|u\|_{\mathcal{B}} \sup_{w\in K}|f_{n}(w)|+C\ve+M\sup_{w\in
K}|f'_{n}(w)|+C\ve+|u(0)f_n(\varphi(0))|,
\end{align*}
where we have used the fact that $u\in {\cal B}$ (see the proof of
Theorem~3.1) and where
\begin{align*}
M=\sup_{z\in D}(1-|z|^2) |u(z)\varphi'(z)|.
\end{align*}
Since $K$ is compact and $\vp\in H(D),$ it follows that,
$\lim_{n\to\infty}\sup_{w\in K}|f_{n}(w)|=0.$ The set $\{\vp(0)\}$
is also compact so that $\lim_{n\to\infty}|u(0) f_{n}(\vp(0))|=0.$
By Cauchy's estimate, if $f_n$ is a sequence which converges on
compacta of $D$ to zero, then the sequence $f'_n$ also converges
on compacta of $D$ to zero as $n\to\infty.$ Employing these facts
and letting $n\to\infty$ in the last inequality, we obtain that
\begin{equation*}
\limsup_{n\to\infty}\|uC_{\varphi}f_n\|_{\mathcal{B}} \leq 2C\ve.
\end{equation*}
Since $\ve$ is an arbitrary positive number it follows that the last limit is equal to
zero. Therefore, $uC_{\varphi}\hbox{\rm :}\ H(p,p,\phi) \rightarrow \mathcal{B}$ is compact.
\end{proof}

Conversely, suppose $uC_{\varphi}\hbox{\rm :}\ H(p,p,\phi) \rightarrow
\mathcal{B}$ is compact. Let $(z_n)_{n\in\NN}$ be a sequence in
$D$ such that $|\varphi(z_n)|\rightarrow 1$ as $n\rightarrow
\infty$. If such a sequence does not exist conditions (15) and
(16) are automatically satisfied. Choose
\begin{equation}
f_n(z)=\frac{(1-|\varphi(z_n)|^{2})^{t+1}}{\phi(|\varphi(z_n)|)
(1-\overline{\varphi(z_n)} z)^{1/p+t+1} },\qquad n\in\NN.
\end{equation}
Then, as above $\sup_{n\in\NN}\|f_n\|_{H(p,p,\phi)}\leq C$ and
$f_n$ converges to 0 uniformly on compact subsets of $D$ as
$n\to\infty.$ Since $uC_{\varphi}$ is compact, we have
\begin{equation*}
\|uC_{\varphi}f_{n}\|_{\mathcal{B}} \rightarrow 0 \quad\mbox{as} \
n\rightarrow\infty.
\end{equation*}
Thus
\begin{align*}
&\hskip -4pc \|uC_{\varphi}f_{n}\|_{\mathcal{B}}\\[.3pc]
&\hskip -4pc \quad\,\geq
\sup_{z \in D} (1-|z|^{2})|(uC_{\varphi}f_{n})'(z)| \\[.4pc]
&\hskip -4pc \quad\,\geq \left|(1/p+t+1) \frac{(1-|z_n|^{2}) |u(z_n)
\overline{\varphi(z_n)} \varphi'(z_n)|}{ \phi(|\varphi(z_n)|)
(1-|\varphi(z_n)|^2)^{1+1/p}} - \frac{(1-|z_n|^2)|u'(z_n)|} {\phi(|\varphi(z_n)|)
(1-|\varphi(z_n)|^2)^{1/p}}\right|.
\end{align*}
Hence, we obtain
\begin{align}
&\lim_{|\varphi(z_n)|\rightarrow 1 } \frac{(1/p+t+1)(1-|z_n|^{2} )|u(z_n)
\overline{\varphi(z_n)} \varphi'(z_n)|}{ \phi(|\varphi(z_n)|) (1-|\varphi(z_n)|^2)^{
1+1/p}}\nonumber\\[.7pc] &\quad\,= \lim_{|\varphi(z_n)|\rightarrow 1 }
\frac{(1-|z_n|^2)|u'(z_n)|} {\phi(|\varphi(z_n)|) (1-|\varphi(z_n)|^2)^{1/p}},
\end{align}
if one of these two limits exists.

Next, let
\begin{equation}
g_n(z)=\frac{(1-|\varphi(z_n)|^{2})^{t+2}}{\phi(|\varphi(z_n)|)
(1-\overline{\varphi(z_n)} z)^{1/p+t+2}
}-\frac{(1-|\varphi(z_n)|^{2})^{t+1}}{\phi(|\varphi(z_n)|)
(1-\overline{\varphi(z_n)} z)^{1/p+t+1} }
\end{equation}
for a sequence $(z_n)_{n\in\NN}$ in $D$ such that
$|\varphi(z_n)|\rightarrow 1$ as $n\to\infty.$ Then,
$(g_n)_{n\in\NN}$ is a bounded sequence in $H(p,p,\phi)$,
$g_n\rightarrow 0$ uniformly on every compact subset of $D$ as
$n\to\infty$, $g_n(\varphi(z_n))=0$ and
\begin{equation*}
g_n'(\varphi(z_n))=\frac{\overline{\varphi(z_n)}}{\phi(|\varphi(z_n)|)(1-|\varphi(z_n)|^2)^{1+1/p}}.
\end{equation*}
Then
\begin{align}
&\frac{(1-|z_n|^2) |u(z_n)\overline{\varphi(z_n)} \varphi'(z_n)|}
 {\phi(|\varphi(z_n)|)(1-|\varphi(z_n)|^2)^{1+\frac{1}{p} } } \leq\|uC_{\varphi}g_n\|_{\mathcal{B}}\to 0
\end{align}
as $n\to\infty.$

From (20) it follows that
\begin{equation*}
\lim_{|\varphi(z)|\rightarrow
1}\frac{(1-|z|^{2})|u(z)\varphi'(z)|}{\phi(|\varphi(z)|)
(1-|\varphi(z)|^{2})^{1+1/p}}=0.
\end{equation*}
Therefore by (18), we have
\begin{equation*}
\lim_{|\varphi(z)|\rightarrow
1}\frac{(1-|z|^2)|u'(z)|} {\phi(|\varphi(z)|)
(1-|\varphi(z)|^2)^{1/p}}=0.
\end{equation*}

From the last two theorems, we can easily obtain the following corollaries:

\begin{coro}\label{abelian}$\left.\right.$\vspace{.5pc}

\noindent Suppose that $\varphi $ is an analytic self-map of the unit disk{\rm ,} $0<
p<\infty$ and that $\phi$ is normal on $[0,1)$. Then{\rm ,} the composition operator
$C_{\varphi}\hbox{\rm :}\ H(p,p,\phi)\rightarrow \mathcal{B}$ is bounded if and only
if the following condition is satisfied{\rm :}
\begin{equation*}
\sup_{z \in D}
\frac{(1-|z|^2)|\varphi'(z)|}{\phi(|\varphi(z)|)(1-|\varphi(z)|^2)^{1+1/p}}<\infty.
\end{equation*}
\end{coro}\pagebreak

\begin{coro}\label{abelian}$\left.\right.$\vspace{.5pc}

\noindent Suppose that $\varphi$ is an analytic self-map of the unit disk{\rm ,} $0<
p<\infty$, that $\phi$ is normal on $[0,1)$ and that $C_{\varphi}\hbox{\rm :}\
H(p,p,\phi)\rightarrow \mathcal{B}$ is bounded. Then{\rm ,} $C_{\varphi}\hbox{\rm :}\
H(p,p,\phi)\rightarrow \mathcal{B}$ is compact if and only if the following condition
is satisfied{\rm :}
\begin{equation*}
\lim_{|\varphi(z)|\rightarrow 1}
\frac{(1-|z|^2)|\varphi'(z)|}{\phi(|\varphi(z)|)
(1-|\varphi(z)|^2)^{1+1/p}}=0.
\end{equation*}
\end{coro}

Since the Bergman space is a special case of $H(p,p,\phi)$, we have the following
corollaries.

\begin{coro}\label{abelian}$\left.\right.$\vspace{.5pc}

\noindent Suppose that $\varphi$ is an analytic self-map of the unit disk{\rm ,} $u\in
H(D)$ and $0<p<\infty$. Then{\rm ,} $uC_{\varphi}\hbox{\rm :}\ A^p\rightarrow
\mathcal{B}$ is bounded if and only if the following  conditions are satisfied{\rm :}
\begin{equation*}\sup_{z\in
D}\frac{(1-|z|^{2})|u'(z)|}{(1-|\varphi(z)|^{2})^{2/p}}<\infty
\qand  \sup_{z\in D}\frac{(1-|z|^{2})|u(z)\varphi'(z)|
}{(1-|\varphi(z)|^{2})^{1+2/p}}<\infty.
\end{equation*}\vspace{.1pc}
\end{coro}

\begin{coro}\label{abelian}$\left.\right.$\vspace{.5pc}

\noindent Suppose that $\varphi $ is an analytic self-map of the unit disk{\rm ,} $u
\in H(D),$ $0<p<\infty$ and that $uC_{\varphi}\hbox{\rm :}\ A^p\rightarrow
\mathcal{B}$ is bounded. Then{\rm ,} $uC_{\varphi}\hbox{\rm :}\ A^p\rightarrow
\mathcal{B}$ is compact if and only if the following  conditions are satisfied{\rm :}
\begin{align*}
\lim_{|\varphi(z)|\rightarrow
1}\frac{(1-|z|^{2})|u'(z)|}{(1-|\varphi(z)|^{2})^{2/p}}=0\qand
\lim_{|\varphi(z)|\rightarrow
1}\frac{(1-|z|^{2})|u(z)\varphi'(z)|}{(1-|\varphi(z)|^{2})^{1+2/p}}=0.
\end{align*}\vspace{.1pc}
\end{coro}

\begin{coro}\label{abelian}$\left.\right.$\vspace{.5pc}

\noindent Suppose that $\varphi $ is an analytic self-map of the unit disk and
$0<p<\infty$. Then{\rm ,} $C_{\varphi}\hbox{\rm :}\ A^p\rightarrow \mathcal{B}$ is
bounded if and only if the following  condition is satisfied{\rm :}
\begin{equation*}
\sup_{z \in D}
\frac{(1-|z|^2)|\varphi'(z)|}{(1-|\varphi(z)|^2)^{1+2/p}}<\infty.
\end{equation*}\vspace{.1pc}
\end{coro}

\begin{coro}\label{abelian}$\left.\right.$\vspace{.5pc}

\noindent Suppose that $\varphi $ is an analytic self-map of the unit disk{\rm ,} $0<
p<\infty$ and that $C_{\varphi}\hbox{\rm :}\ A^p\rightarrow \mathcal{B}$ is bounded.
Then{\rm ,} $C_{\varphi}\hbox{\rm :}\ A^p\rightarrow \mathcal{B}$ is compact if and
only if the following  condition is\break satisfied{\rm :}
\begin{equation*}
\lim_{|\varphi(z)|\rightarrow 1}
\frac{(1-|z|^2)|\varphi'(z)|}{(1-|\varphi(z)|^2)^{1+2/p}}=0.
\end{equation*}
\end{coro}

\section{The boundedness and compactness of the operator \pmb{$uC_{\varphi}\hbox{\rm :}\ H(p,p,\phi)\rightarrow
\mathcal{B}_0$}}

Next we characterize the boundedness and compactness of the
weighted composition operators $uC_{\varphi}\hbox{\rm :}\
H(p,p,\phi) \rightarrow \mathcal{B}_0$. For this purpose, we need
the following lemmas. The first lemma can be found in \cite{mm}.
\pagebreak

\setcounter{theore}{0}
\begin{lem} A closed set K in $\mathcal{B}_0$ is
compact if and only if it is bounded and satisfies
\begin{equation}
\lim_{|z|\rightarrow 1}\sup_{f \in K}(1-|z|^2)|f'(z)|=0.
\end{equation}
\end{lem}

\begin{lem} Suppose that $\varphi$ is an
analytic self-map of the unit disk{\rm ,} $u \in H(D)$, $0< p<\infty$ and that $\phi$
is normal on $[0,1)$. Then{\rm ,}
\begin{equation}\lim_{|z|\rightarrow
1}\frac{(1-|z|^{2})|u'(z)|}{\phi(|\varphi(z)|)(1-|\varphi(z)|^{2})^{1/p}}=0
\end{equation}
if and only if
\begin{equation}\lim_{|\varphi(z)|\rightarrow
1}\frac{(1-|z|^2)|u'(z)|} {\phi(|\varphi(z)|)
(1-|\varphi(z)|^2)^{1/p}}=0
\end{equation}
and
\begin{align}u\in \mathcal{B}_0.
\end{align}
\end{lem}

\begin{proof}
Suppose that (22) holds. Then
\begin{equation*}
(1-|z|^{2})|u'(z)|\leq
\frac{C(1-|z|^{2})|u'(z)|}{\phi(|\varphi(z)|)(1-|\varphi(z)|^{2})^{1/p}}\rightarrow
0
\end{equation*}
as $|z| \rightarrow 1$.

If $|\varphi(z)|\rightarrow 1$, then $|z| \rightarrow 1$, from
which it follows that
\begin{equation*}
\lim_{|\varphi(z)|\rightarrow 1}\frac{(1-|z|^2)|u'(z)|}
{\phi(|\varphi(z)|) (1-|\varphi(z)|^2)^{1/p}}=0.
\end{equation*}

Conversely, suppose that (23) and (24) hold. By (23), for every
$\epsilon>0$, there exists $r\in (0,1)$,
\begin{equation*}
\frac{(1-|z|^2)|u'(z)|} {\phi(|\varphi(z)|)
(1-|\varphi(z)|^2)^{1/p}}<\epsilon
\end{equation*}
when $r<|\varphi(z)|<1$. By (24), there exists $\sigma \in (0,1)$,
\begin{equation*}
(1-|z|^{2})|u'(z)|\leq \epsilon(1-r^2)^{1/p}\phi(r)
\end{equation*}
when $\sigma<|z|<1.$

Therefore, when $\sigma<|z|<1$ and $r<|\varphi(z)|<1$, we have
that
\begin{equation}
\frac{(1-|z|^2)|u'(z)|} {\phi(|\varphi(z)|)
(1-|\varphi(z)|^2)^{1/p}}<\epsilon.
\end{equation}
If $ |\varphi(z)|\leq r$ and $\sigma<|z|<1,$ then since $\phi$ is normal, we obtain
\begin{equation}\frac{(1-|z|^2)|u'(z)|} {\phi(|\varphi(z)|)
(1-|\varphi(z)|^2)^{1/p}}<\frac{(1-r)^s(1-|z|^2)|u'(z)|} {\phi(r)
(1-|\varphi(z)|^2)^{1/p+s}}<\epsilon.
\end{equation}
Combining (25) with (26), we obtain the desired result.
\end{proof}

Similarly to the proof of the above lemma, we have the following.\pagebreak

\begin{lem}
Suppose that $\varphi $ is an analytic self-map of the unit disk{\rm ,} $u \in
H(D)${\rm ,} $0< p<\infty$ and that $\phi$ is normal on $[0,1)$. Then{\rm ,}
\begin{equation}\lim_{|z|\rightarrow
1}\frac{(1-|z|^{2})|u(z)\varphi'(z)|}{\phi(|\varphi(z)|)(1-|\varphi(z)|^{2})^{1+1/p}}=0
\end{equation}
if and only if
\begin{equation}\lim_{|\varphi(z)|\rightarrow
1}\frac{(1-|z|^{2})|u(z)\varphi'(z)|}{\phi(|\varphi(z)|)
(1-|\varphi(z)|^{2})^{1+1/p}}=0
\end{equation}
and
\begin{equation}
\lim_{|z|\to 1}(1-|z|)^2|u(z)\vp'(z)|=0.
\end{equation}
\end{lem}

\begin{theor}[\!]
Suppose that $\varphi$ is an analytic self-map of the unit disk{\rm ,} $u \in H(D),$
$0< p<\infty$ and that $\phi$ is normal on $[0,1)$. Then{\rm ,} $uC_{\varphi}\hbox{\rm
:}\ H(p,p,\phi) \rightarrow \mathcal{B}_0$ is bounded if and only if
$uC_{\varphi}\hbox{\rm :}\ H(p,p,\phi) \rightarrow \mathcal{B}$ is bounded{\rm ,}
$u\in \mathcal{B}_0$ and
\begin{equation*}
\lim_{|z|\to 1}(1-|z|^2)|u(z)\vp'(z)|=0.
\end{equation*}
\end{theor}

\begin{proof}
First assume that $uC_{\varphi}\hbox{\rm :}\
H(p,p,\phi)\rightarrow \mathcal{B}_0$ is bounded. Then, it is
clear that $uC_{\varphi}\hbox{\rm :}\ H(p,p,\phi)\rightarrow \mathcal{B}$ is
bounded. Taking the functions $f(z)=1$ and $f(z)=z,$ we obtain
that $u\in \mathcal{B}_0$ and $\lim_{|z|\to
1}(1-|z|^2)|u(z)\vp'(z)|=0.$

Conversely, assume that $uC_{\varphi}\hbox{\rm :}\ H(p,p,\phi) \rightarrow
\mathcal{B}$ is bounded, $u \in \mathcal{B}_0$ and $\lim_{|z|\to
1}(1-|z|^2)|u(z)\vp'(z)|=0$. Then, for each polynomial $p(z)$, we
have that
\begin{align*}
&(1-|z|^2) |(uC_\varphi p)'(z)|\\[.4pc]
&\quad\,\leq  (1-|z|^2) |u'(z)| | p( \varphi (z))|+
(1-|z|^2)|u(z)\varphi'(z)p'(\varphi(z))|,
\end{align*}
from which it follows that $uC_\varphi p\in \mathcal{B}_0$. Since
the set of all polynomials is dense in $H(p,p,\phi),$ we have that
for every $f\in H(p,p,\phi)$ there is a sequence of polynomials
$(p_n)_{n\in\NN}$ such that $\| f-p_n\|_{H(p,p,\phi) }\rightarrow
0$, as $n\rightarrow\infty.$ Hence
\begin{equation*}\|uC_\vp f-uC_\vp p_n\|_{\cal B}\leq \|uC_\vp\|_{H(p,p,\phi)\to\mathcal{B}}\|f-p_n\|_{H(p,p,\phi)}\to 0\end{equation*}
as $n\to\infty,$ since the operator $uC_{\varphi}\hbox{\rm :}\ H(p,p,\phi)
\rightarrow \mathcal{B}$ is bounded. Since $\mathcal{B}_0$ is a
closed subset of $\mathcal{B}$, we obtain
\begin{equation*}
uC_\varphi (H(p,p,\phi))\subset\mathcal{B}_0.
\end{equation*}
Therefore $ uC_\varphi\hbox{\rm :}\ H(p,p,\phi) \rightarrow \mathcal{B}_0$ is bounded.
\end{proof}

\begin{theor}[\!]
Suppose that $\varphi$ is an analytic self-map of the unit disk{\rm ,} $u \in H(D),$
$0< p<\infty$ and that $\phi$ is normal on $[0,1)$. Then{\rm ,} $uC_{\varphi}\hbox{\rm
:}\ H(p,p,\phi)\rightarrow \mathcal{B}_0$ is compact if and\break only if
\begin{align*}
\hskip -4pc \lim_{|z|\rightarrow
1}\frac{(1-|z|^{2})|u'(z)|}{\phi(|\varphi(z)|)(1-|\varphi(z)|^{2})^{1/p}}=0\qand
   \lim_{|z|\rightarrow
1}\frac{(1-|z|^{2})|u(z)\varphi'(z)|}{\phi(|\varphi(z)|)(1-|\varphi(z)|^{2})^{1+1/p}}=0.
\end{align*}
\end{theor}

\begin{proof}
First, we assume that $uC_{\varphi}\hbox{\rm :}\
H(p,p,\phi)\rightarrow \mathcal{B}_0$ is compact. Taking $f(z)\equiv
1$ we obtain that
\begin{equation}
u\in \mathcal{B}_0.
\end{equation}
From this, taking $f(z)=z,$ and using the boundedness of $uC_{\varphi}\hbox{\rm :}\
H(p,p,\phi)\rightarrow \mathcal{B}_0$ it follows that
\begin{equation}
\lim_{|z|\to 1}(1-|z|)^2|u(z)\vp'(z)|=0.
\end{equation}
Hence, if $\|\vp\|_\infty<1,$ from (30) and (31), we obtain that
\begin{equation*}
\lim_{|z|\rightarrow
1}\frac{(1-|z|^{2})|u'(z)|}{\phi(|\varphi(z)|)(1-|\varphi(z)|^{2})^{1/p}}\leq
C\lim_{|z|\rightarrow
1}\frac{(1-|z|^{2})|u'(z)|}{\phi(\|\vp\|_\infty)(1-\|\vp\|_\infty^{2})^{1/p}}=0
\end{equation*}
and
\begin{align*}
&\lim_{|z|\rightarrow
1}\frac{(1-|z|^{2})|u(z)\varphi'(z)|}{\phi(|\varphi(z)|)(1-|\varphi(z)|^{2})^{1+1/p}}\\[.45pc]
&\quad\,\leq C\lim_{|z|\rightarrow
1}\frac{(1-|z|^{2})|u(z)\varphi'(z)|}{\phi(\|\vp\|_\infty)(1-\|\vp\|_\infty^{2})^{1+1/p}}=0
\end{align*}
from which the conditions in (22) and (27) follow.

Hence, assume that $\|\vp\|_\infty=1.$ Let $(\vp(z_n))_{n\in\NN}$ be a sequence such
that $\lim_{n\to\infty}|\vp(z_n)|$ $=1$. If necessary we can take a subsequence of
$(\vp(z_n))_{n\in\NN}$ (we use the same notation $(\vp(z_n))_{n\in\NN}$). Set
\begin{equation*}
f_n(z)=\frac{(1-|\varphi(z_n)|^{2})^{t+1}}{\phi(|\varphi(z_n)|)
(1-\overline{\varphi(z_n)} z)^{1/p+t+1} },\qquad n\in\NN
\end{equation*}
and
\begin{equation*}
g_n(z)=\frac{(1-|\varphi(z_n)|^{2})^{t+2}}{\phi(|\varphi(z_n)|)
(1-\overline{\varphi(z_n)} z)^{1/p+t+2}
}-\frac{(1-|\varphi(z_n)|^{2})^{t+1}}{\phi(|\varphi(z_n)|) (1-\overline{\varphi(z_n)}
z)^{1/p+t+1}}.
\end{equation*}
By the proof of Theorem~3.2 we know that
\begin{equation}\lim_{|\varphi(z)|\rightarrow
1}\frac{(1-|z|^{2})|u(z)\varphi'(z)|}{\phi(|\varphi(z)|)
(1-|\varphi(z)|^{2})^{1+1/p}}=0
\end{equation}
and
\begin{equation}
\lim_{|\varphi(z)|\rightarrow
1}\frac{(1-|z|^2)|u'(z)|} {\phi(|\varphi(z)|)
(1-|\varphi(z)|^2)^{1/p}}=0.
\end{equation}
Applying (30), (31), (32) and (33) with Lemmas~4.2 and 4.3 gives
the desired result.

Conversely, from (7) we have that
\begin{align}
&\hskip -4pc (1-|z|^{2})|(uC_{\varphi}f)'(z)|\nonumber\\[.45pc]
&\hskip -4pc \quad\,\leq C\left( \frac{(1-|z|^2)|u'(z)|}{\phi(|\varphi(z)|)
(1-|\varphi(z)|^2)^{1/p} }+ \frac{(1-|z|^2) |u(z)\varphi'(z)|
}{\phi(|\varphi(z)|)(1-|\varphi(z)|^2)^{1+1/p} } \right)
 \|f\|_{H(p,p,\phi)}.\nonumber
\end{align}
Taking the supremum in this inequality over all $f\in H(p,p,\phi)$
such that $\|f\|_{H(p,p,\phi)}\leq 1,$ then letting $|z|\to 1$, we
obtain that
\begin{equation*}
\lim_{|z|\rightarrow 1}\sup_{\|f\|_{H(p,p,\phi)}\leq
1}(1-|z|^2)|(u C_\vp(f))'(z)|=0,
\end{equation*}
from which by Lemma~4.1 we obtain that the operator $uC_{\varphi}\hbox{\rm :}\
H(p,p,\phi)\rightarrow \mathcal{B}_0$ is compact.
\end{proof}

From Theorems~4.4 and 4.5, we obtain the following corollaries:

\begin{coro}\label{abelian}$\left.\right.$\vspace{.5pc}

\noindent Suppose that $\varphi$ is an analytic self-map of the unit disk{\rm ,} $0<
p<\infty$ and that $\phi$ is normal on $[0,1)$. Then{\rm ,} the following statements
hold.

\begin{enumerate}
\renewcommand\labelenumi{{\rm (\roman{enumi})}}
\leftskip .2pc
\item $C_{\varphi}\hbox{\rm :}\ H(p,p,\phi) \rightarrow \mathcal{B}_0$ is bounded
if and only if $C_{\varphi}\hbox{\rm :}\ H(p,p,\phi) \rightarrow \mathcal{B}$ is
bounded and $\varphi \in \mathcal{B}_0$.

\item $C_{\varphi}\hbox{\rm :}\ H(p,p,\phi) \rightarrow \mathcal{B}_0$ is compact if
and only if
\begin{equation*}
\hskip -1.25pc \lim_{|z|\rightarrow
1}\frac{(1-|z|^{2})|\varphi'(z)|}{\phi(|\varphi(z)|) (1-|\varphi(z)|^{2})^{1+1/p}}=0.
\end{equation*}
\end{enumerate}
\end{coro}

\begin{coro}\label{abelian}$\left.\right.$\vspace{.5pc}

\noindent Suppose that $\varphi$ is an analytic self-map of the unit disk{\rm ,} $u
\in H(D)$ and $0< p<\infty$. Then{\rm ,} the following statements hold.

\begin{enumerate}
\renewcommand\labelenumi{{\rm (\roman{enumi})}}
\leftskip .2pc
\item $uC_{\varphi}\hbox{\rm :}\ A^p \rightarrow \mathcal{B}_0$ is bounded if and only
if $ uC_{\varphi}\hbox{\rm :}\ A^p \rightarrow \mathcal{B}$ is bounded{\rm ,} $u \in
\mathcal{B}_0$ and
\begin{equation*}
\hskip -1.25pc \lim_{|z|\rightarrow1} (1-|z|^{2})|u(z)\varphi'(z)|=0.
\end{equation*}

\item $uC_{\varphi}\hbox{\rm :}\ A^p \rightarrow \mathcal{B}_0$ is compact if and only
if
\begin{align*}
\hskip -1.25pc \lim_{|z|\rightarrow
1}\frac{(1-|z|^{2})|u'(z)|}{(1-|\varphi(z)|^{2})^{2/p}}=0 \ \ \ \,\hbox{and}\, \ \ \
\lim_{|z|\rightarrow
1}\frac{(1-|z|^{2})}{(1-|\varphi(z)|^{2})^{1+2/p}}|u(z)\varphi'(z)|=0.
\end{align*}
\end{enumerate}
\end{coro}

\begin{coro}\label{abelian}$\left.\right.$\vspace{.5pc}

\noindent Suppose that $\varphi$ is an analytic self-map of the unit disk and
$0<p<\infty$. Then{\rm ,} the following statements hold.

\begin{enumerate}
\renewcommand\labelenumi{{\rm (\roman{enumi})}}
\leftskip .2pc
\item $C_{\varphi}\hbox{\rm :}\ A^p \rightarrow \mathcal{B}_0$ is bounded if and
only if $C_{\varphi}\hbox{\rm :}\ A^p \rightarrow \mathcal{B}$ is bounded and
 $\varphi \in \mathcal{B}_0$.

\item $ C_{\varphi}\hbox{\rm :}\ A^p \rightarrow \mathcal{B}_0$ is compact if and only
if
\begin{equation*}
\hskip -1.25pc \lim_{|z|\rightarrow
1}\frac{(1-|z|^{2})}{(1-|\varphi(z)|^{2})^{1+2/p}}|\varphi'(z)|=0.\end{equation*}
\end{enumerate}
\end{coro}

\section*{Acknowledgement}

The first author of this paper is supported in part by the NNSF of
China (No.~10671115), grants from specialized research fund for
the doctoral program of higher education (No.~20060560002) and NSF
of Guangdong Province (No.~06105648).

\end{document}